\documentclass[11pt,a4paper]{article}
\usepackage{pifont}
\usepackage{CJK,epsfig,amsfonts,amsgen,amsmath,amstext,amsbsy,amsopn,amsthm
}
\usepackage{mathrsfs}
\usepackage{float}
\usepackage{epsfig,latexsym,amsfonts,amsmath,amssymb}
\usepackage{color}
\usepackage{calc}
\usepackage{tikz}

\usepackage{comment}
\usepackage{geometry}
\usepackage{calc}
\usepackage{tikz}
\usetikzlibrary{decorations.markings}
\tikzstyle{vertex}=[circle, draw, inner sep=0pt, minimum size=6pt]

\newtheorem{theorem}[subsection]{Theorem}
\newtheorem{lemma}[subsection]{Lemma}

\newtheorem{remark}[subsection]{Remark}

\newtheorem{case}[]{Case}

\usepackage{indentfirst} 
\newcounter{Hcase}

\begin{document}
\title{Largest family without a pair of posets \\on consecutive levels of the Boolean lattice}
\author{\small Gyula O.H. Katona${}^{1}$ and Jimeng Xiao${}^{2,3}$
 \\[2mm]
\small ${}^1$ Alfr\'ed R\'enyi Institute of Mathematics, Budapest\\
\small ${}^2$School of Mathematics and Statistics, \\
\small Northwestern Polytechnical University, Xi'an, P.R.~China\\
\small ${}^3$ Xi'an-Budapest Joint Research Center for Combinatorics,\\ \small Northwestern Polytechnical University, Xi'an, P.R. China}

\baselineskip 20pt

\date{}
\maketitle
 \vspace{4mm}

\begin{abstract}
 
Suppose $k \ge 2$ is an integer. 
Let $Y_k$ be the poset with elements $x_1, x_2, y_1, y_2, \ldots, y_{k-1}$ such that $y_1 < y_2 < \cdots < y_{k-1} < x_1, x_2$ and let $Y_k'$ be the same poset but all relations reversed. 
We say that a family of subsets of $[n]$ contains a copy of $Y_k$ on consecutive levels if it contains $k+1$ subsets $F_1, F_2, G_1, G_2, \ldots, G_{k-1}$ such that  $G_1\subset G_2 \subset \cdots \subset G_{k-1} \subset F_1, F_2$ and $|F_1| = |F_2| = |G_{k-1}|+1 =|G_{k-2}|+ 2= \cdots = |G_{1}|+k-1$. 
If both $Y_k$ and $Y'_k$ on consecutive levels are forbidden, the size of the largest such family is denoted by $\mathrm{La}_{\mathrm{c}}(n, Y_k, Y'_k)$. In this paper, we will determine the exact value of $\mathrm{La}_{\mathrm{c}}(n, Y_k, Y'_k)$.

\vskip 0.1in \noindent%
 \\[7pt]
\end{abstract}

\newcommand{\lr}[1]{\langle #1\rangle}

\parindent 17pt
\baselineskip 17pt

\section{Introduction}

Given two partially ordered sets (posets) $P$ and $Q$, we say that $P$ is a \emph{subposet} of $Q$ if there exists an injection $\phi:P \to Q$ such that $x \le_P y$ implies $\phi(x)~\le_Q~\phi(y)$.   
Viewing collections of sets as posets under the inclusion relation, we have the following extremal functions, first introduced by Katona and Tarj\'an \cite{kt}. 
For any collection of finite posets $\mathcal{P}$, let $\mathrm{La}(n, \mathcal{P})$ be the maximum size of a family of subsets of $[n]=\{1,2,\ldots,n\}$ which does not contain any $P \in \mathcal{P}$ as a subposet.

This type of problems was first studied by Sperner \cite{sp}.  
\begin{theorem}[Sperner \cite{sp}]
Let $\mathcal{F}$ be a family of subsets of $[n]$ without inclusion relation between any two of the subsets. 
Then
\[
|\mathcal{F}| \le \binom{n}{\lfloor\frac{n}{2}\rfloor}.
\]
\end{theorem}

A chain of length $k$ is a poset with elements $x_1, x_2, \ldots, x_k$ such that $x_1< x_2<  \cdots< x_k$. Let $\sum(n,k)$ be the sum of the $k$ largest binomial coefficients of the form $\binom{n}{i}$. Then Sperner theorem can be extended as follows:
\begin{theorem}[Erd\H{o}s \cite{e}]
Let $\mathcal{F}$ be a family of subsets of $[n]$ without a chain of length $k$. 
Then
\[
|\mathcal{F}| \le \sum(n,k-1).
\]
\end{theorem}

Before stating the next result, we need the following notation. 
Let $2 \le k \le n$ and $0\le r \le n$ be two integers. 
The following lacunary sum of binomial coefficients was first introduced by Ramus \cite{r} in 1834.
\[
 S(n,k,r) = \sum_{\substack{i=0\\ i \equiv r \text{~mod~} k}}^{n}\binom{n}{i}.
\]
Clearly, if $r_1 \equiv r_2$ (mod $k$), then $S(n,k,r_1) = S(n,k,r_2)$. So there are $k$ distinct such sums $S(n,k,0), S(n,k,1), \ldots, S(n,k,k-1)$.

The first author \cite{k} published the next two theorems which are analogs of the two theorems above.
\begin{theorem}[Katona \cite{k}]
Let $\mathcal{F}$ be a family of subsets of $[n]$ such that no two of the subsets $F_i, F_j$ satisfy $F_i \subset F_j$ and $|F_j| - |F_i| < k$. Then
\[
|\mathcal{F}| \le \max\{S(n,k,r) \mid r \in \{0, \ldots, k-1\}\}.
\]
\end{theorem}
\begin{theorem}[Katona \cite{k}]
Let $h \le k$ be an integer, and let $\mathcal{F}$ be a family of subsets of $[n]$ such that no $h+1$ of the subsets $F_{i_{1}}, \ldots, F_{i_{h+1}}$ satisfy $F_{i_1} \subset \cdots \subset F_{i_{h+1}}$ and $|F_{i_{h+1}}| - |F_{i_1}| < k$. 
Then
\[
|\mathcal{F}| \le \max\{S(n,k,r_1) + S(n,k,r_2) + \cdots + S(n,k,r_h)  \mid r_1 \not= r_2 \not= \cdots  \not= r_h \in \{0, \ldots, k-1\}\}.
\]
\end{theorem}

Another type of generalizations of Sperner Theorem is to determine the largest size of  a family of  subsets of $[n]$ without a  copy of poset $Y_k$ and $Y'_k$ defined below. 
Let $Y_k$ be the poset with elements $x_1, x_2, y_1, y_2, \ldots, y_{k-1}$ such that $y_1 < y_2 < \cdots < y_{k-1} < x_1, x_2$, and let $Y'_k$ be the same poset but all relations reversed. 
Katona and Tarj\'an \cite{kt} gave the following result for $k = 2$. (In their paper, posets $Y_2$ and $Y'_2$ are denoted by $\mathrm{V}$ (the cheery poset) and $\Lambda$ (the fork poset) respectively.)
\begin{theorem}[Katona and Tarj\'an \cite{kt}]
\[
\mathrm{La}(n, Y_2, Y'_2) = 2\binom{n-1}{\lfloor\frac{n-1}{2}\rfloor}.
\]
\end{theorem}

De Bonis, Katona and Swanepoel \cite{dk} studied the case $k=3$. We remark that they actually proved a result for a poset so called butterfly, but their proof implies the following theorem.
\begin{theorem}[De Bonis, Katona and Swanepoel \cite{dk}]
\[
\mathrm{La}(n, Y_3, Y'_3) = \sum(n,2).
\]
\end{theorem}

For the case $k \ge 4$, Methuku and Tompkins \cite{mt} got the next theorem.
\begin{theorem}[Methuku and Tompkins \cite{mt}]
For $n \ge k \ge 4$,
\[
\mathrm{La}(n, Y_k, Y'_k) = \sum(n,k-1).
\]
\end{theorem}

For other results related to the La function, see \cite{bn,gm,mm,tw}. 
Now, we consider the following problem, which is a combination of the two types of problems above. 
We say that a family of subsets of $[n]$ contains a copy of $Y_k$ on consecutive levels if it contains $k+1$ subsets $F_1, F_2, G_1, G_2, \ldots, G_{k-1}$ such that  $G_1\subset G_2 \subset \cdots \subset G_{k-1} \subset F_1, F_2$ and $|F_1| = |F_2| = |G_{k-1}|+1 =|G_{k-2}|+ 2= \cdots = |G_{1}|+k-1$.
We denote by $\mathrm{La}_{\mathrm{c}}(n, Y_k, Y'_k)$ the largest size of a family of subsets of $[n]$ without both $Y_k$ and $Y'_k$ on consecutive levels.
We note that this problem is mentioned in \cite{gmn} in the language of the oriented hypercube. The authors in \cite{gmn} gave  an asymptotic formula of the size of largest family without tree posets in consecutive levels. 
For exact result, they proved 
\[
\mathrm{La}_{\mathrm{c}}(n, Y_2, Y'_2) =2^{n-1}.
\]
Let $m = \lceil(n-k)/2\rceil$ in the rest of this paper. For $k \ge 3$, we have the following theorem.
\begin{theorem}\label{thm}
Let $n \ge k \ge 3$, then
\[
\mathrm{La}_{\mathrm{c}}(n, Y_k, Y'_k) =2^{n} - S(n,k,m).
\]
\end{theorem}
\begin{remark}\label{rem_lm}
Loehr and Michael \cite{lm} showed that 
\[
S(n,k,m) = \min_{r: \ 0\le r \le k-1}S(n,k,r).
\]
So our result implies that the trivial construction is the best. Here, a trivial construction consists of all subsets except the ones with size $s \equiv m$ (mod $k$).
\end{remark}

The rest of the paper is organized as follows. In the next section, we present some preliminary results. In Section 3, we will prove our main theorem (Theorem \ref{lem_thm}), which implies Theorem \ref{thm}.

\section{Preliminary results} 

A cyclic permutation $\sigma$ of $[n]$ is a cyclic ordering $a_1, a_2, \ldots, a_n, a_1$, where $a_i \in [n]$ for $i = 1,2, \ldots, n$. 
Let $\mathcal{F}$ be a family without $Y_k$ and $Y'_k$ on consecutive levels. We say a set $F \in \sigma$ if $F$ is an interval along the cyclic permutation $\sigma$. 

Now, we double count the sum $S$ of  $\phi_{F}$  over all cyclic permutations $\sigma$ and $F \in \mathcal{F}$ such that $F \in \sigma$,
where
\begin{align*}
\phi_{F} = 
\begin{cases}
\binom{n}{|F|},&\mbox{if $F \not= \emptyset$ and $F \not= [n]$;}\\
n, &\mbox{if $F = \emptyset$ or $F = [n]$.}
\end{cases}
\end{align*}

\begin{figure}

\begin{tikzpicture}

\draw [fill] (-5,1) circle [radius = 0.1];

\draw [fill] (-3,1) circle [radius = 0.1];

\draw [fill] (-1,1) circle [radius = 0.1];

\draw [fill] (1,1) circle [radius = 0.1];

\draw [fill] (3,1) circle [radius = 0.1];

\draw [fill] (5,1) circle [radius = 0.1];

\draw [fill] (-5,-1) circle [radius = 0.1];

\draw [fill] (-3,-1) circle [radius = 0.1];

\draw [fill] (-1,-1) circle [radius = 0.1];

\draw [fill] (1,-1) circle [radius = 0.1];

\draw [fill] (3,-1) circle [radius = 0.1];

\draw [fill] (5,-1) circle [radius = 0.1];

\draw [fill] (-5,2) circle [radius = 0.1];

\draw [fill] (-3,2) circle [radius = 0.1];

\draw [fill] (-1,2) circle [radius = 0.1];

\draw [fill] (1,2) circle [radius = 0.1];

\draw [fill] (3,2) circle [radius = 0.1];

\draw [fill] (5,2) circle [radius = 0.1];

\draw [fill] (-5,-2) circle [radius = 0.1];

\draw [fill] (-3,-2) circle [radius = 0.1];

\draw [fill] (-1,-2) circle [radius = 0.1];

\draw [fill] (1,-2) circle [radius = 0.1];

\draw [fill] (3,-2) circle [radius = 0.1];

\draw [fill] (5,-2) circle [radius = 0.1];

\draw [fill] (5,0) circle [radius = 0.1];

\draw [fill] (3,0) circle [radius = 0.1];

\draw [fill] (1,0) circle [radius = 0.1];

\draw [fill] (-1,0) circle [radius = 0.1];

\draw [fill] (-3,0) circle [radius = 0.1];

\draw [fill] (-5,0) circle [radius = 0.1];

\draw [fill] (0,3) circle [radius = 0.1];

\draw [fill] (0,-3) circle [radius = 0.1];

\draw [fill] (0,-3) circle [radius = 0.1];

\node[left] at(-5,-2){\tiny{$\{1\}$}};
\node[right] at(-2.8,-2){\tiny{$\{2\}$}};
\node[right] at(-1,-2){\tiny{$\{3\}$}};

\node[right] at(1,-2){\tiny{$\{4\}$}};
\node[right] at(3,-2){\tiny{$\{5\}$}};
\node[right] at(5,-2){\tiny{$\{6\}$}};

\node[left] at(-5,-1){\tiny{$\{1,2\}$}};
\node[right] at(-2.9,-1){\tiny{$\{2,3\}$}};
\node[right] at(-0.9,-1){\tiny{$\{3,4\}$}};

\node[right] at(1.1,-1){\tiny{$\{4,5\}$}};
\node[right] at(3.05,-0.95){\tiny{$\{5,6\}$}};
\node[right] at(5,-1){\tiny{$\{6,1\}$}};

\node[left] at(-5,0){\tiny{$\{1,2,3\}$}};
\node[right] at(-2.9,0){\tiny{$\{2,3,4\}$}};
\node[right] at(-0.9,0){\tiny{$\{3,4,5\}$}};

\node[right] at(1.1,0){\tiny{$\{4,5,6\}$}};
\node[right] at(3.05,0.05){\tiny{$\{5,6,1\}$}};
\node[right] at(5,0){\tiny{$\{6,1,2\}$}};

\node[left] at(-5,1){\tiny{$\{1,2,3,4\}$}};
\node[right] at(-2.9,1){\tiny{$\{2,3,4,5\}$}};
\node[right] at(-0.9,1){\tiny{$\{3,4,5,6\}$}};

\node[right] at(1.1,1){\tiny{$\{4,5,6,1\}$}};
\node[right] at(3.05,1.05){\tiny{$\{5,6,1,2\}$}};
\node[right] at(5,1){\tiny{$\{6,1,2,3\}$}};

\node[left] at(-5,2){\tiny{$\{1,2,3,4,5\}$}};
\node[right] at(-2.8,1.95){\tiny{$\{2,3,4,5,6\}$}};
\node[right] at(-0.95,2.05){\tiny{$\{3,4,5,6,1\}$}};

\node[right] at(1,2.04){\tiny{$\{4,5,6,1,2\}$}};
\node[right] at(3,2){\tiny{$\{5,6,1,2,3\}$}};
\node[right] at(5,2){\tiny{$\{6,1,2,3,4\}$}};

\node[below] at(0,-3){$\emptyset$};

\node[above] at(0,3){$[6]$};

%

%
\draw[][very thick] (-3,1) -- (1,-1); 
\draw[][very thick] (-5,1) -- (-1,-1); 
\draw[][very thick] (-1,1) -- (3,-1); 

\draw[][very thick] (-5,0) -- (-3,-1); 
\draw[][very thick] (1,1) -- (5,-1); 
\draw[][very thick] (3,1) -- (5,0); 

\draw[][very thick] (-5,0) -- (5,1);
\draw[][very thick] (-5,-1) -- (5,0);
\draw[-][very thick] (-5,1) -- (-5,-1); 
\draw[-][very thick] (-3,1) -- (-3,-1); 
\draw[-][very thick] (-1,1) -- (-1,-1); 
\draw[-][very thick] (1,1) -- (1,-1); 
\draw[-][very thick] (3,1) -- (3,-1); 
\draw[-][very thick] (5,1) -- (5,-1); 
%

\draw[-][very thick] (-5,2) -- (-5,1); 
\draw[-][very thick] (-3,2) -- (-3,1); 
\draw[-][very thick] (-1,2) -- (-1,1); 
\draw[-][very thick] (1,2) -- (1,1); 
\draw[-][very thick] (3,2) -- (3,1); 
\draw[-][very thick] (5,2) -- (5,1); 
\draw[-][very thick] (-5,1) -- (5,2); 
\draw[-][very thick] (-5,2) -- (-3,1); 
\draw[-][very thick] (-3,2) -- (-1,1); 
\draw[-][very thick] (-1,2) -- (1,1); 
\draw[-][very thick] (1,2) -- (3,1); 
\draw[-][very thick] (3,2) -- (5,1); 
\draw[-][very thick] (-5,-2) -- (-5,-1); 
\draw[-][very thick] (-3,-2) -- (-3,-1); 
\draw[-][very thick] (-1,-2) -- (-1,-1); 
\draw[-][very thick] (1,-2) -- (1,-1); 
\draw[-][very thick] (3,-2) -- (3,-1); 
\draw[-][very thick] (5,-2) -- (5,-1); 
\draw[-][very thick] (-5,-2) -- (5,-1); 
\draw[-][very thick] (-5,-1) -- (-3,-2); 
\draw[-][very thick] (-3,-1) -- (-1,-2); 
\draw[-][very thick] (-1,-1) -- (1,-2); 
\draw[-][very thick] (1,-1) -- (3,-2); 
\draw[-][very thick] (3,-1) -- (5,-2); 

\draw[-][very thick] (-5,2) -- (0,3); 

\draw[-][very thick] (-3,2) -- (0,3); 

\draw[-][very thick] (-1,2) -- (0,3); 

\draw[-][very thick] (1,2) -- (0,3); 

\draw[-][very thick] (3,2) -- (0,3); 

\draw[-][very thick] (5,2) -- (0,3);

\draw[-][very thick] (-5,-2) -- (0,-3); 

\draw[-][very thick] (-3,-2) -- (0,-3); 

\draw[-][very thick] (-1,-2) -- (0,-3); 

\draw[-][very thick] (1,-2) -- (0,-3); 

\draw[-][very thick] (3,-2) -- (0,-3); 

\draw[-][very thick] (5,-2) -- (0,-3);

\end{tikzpicture}
\caption{$\mathcal{I}(n)^{\sigma}$ for $n = 6$ and $\sigma = 1$, $2$, $3$, $4$, $5$, $6$, $1$ is shown by bold vertices, and two vertices from consecutive levels are connected by an edge if they have inclusion relation.}\label{fig_eg}
\end{figure}
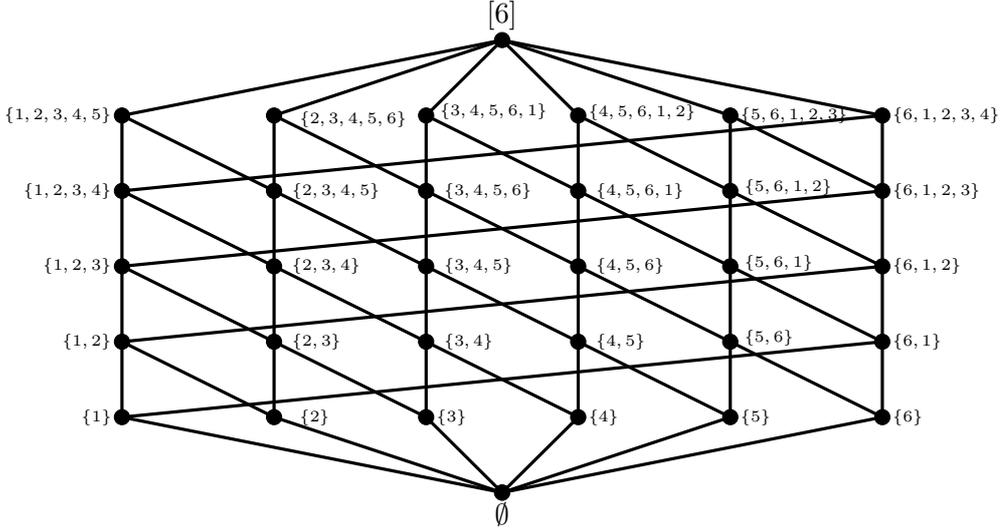

For any $F \in \left(\mathcal{F}\setminus \{\emptyset,[n]\}\right)$, we have $|F|!(n-|F|)!$ cyclic permutations $\sigma$ satisfying $F \in \sigma$. 
If $\emptyset \in \mathcal{F}$ or $[n] \in \mathcal{F}$, all cyclic permutations satisfy the above condition, and the number of cyclic permutations is $(n-1)!$. So
\begin{align*}
  S= &|\{\emptyset,[n]\} \cap \mathcal{F}|\cdot n\cdot(n-1)! + \sum_{F\in  \left(\mathcal{F}\setminus \{\emptyset,[n]\}\right)} |F|!(n-|F|)!\binom{n}{|F|}\\
= &\left(|\{\emptyset,[n]\} \cap \mathcal{F}| + \sum_{F\in  \left(\mathcal{F}\setminus \{\emptyset,[n]\}\right)}1\right) \cdot n! = |\mathcal{F}| \cdot n!.
\end{align*}
On the other hand, for any cyclic permutation $\sigma$,  let $\mathcal{I}(n)^{\sigma}$ be the family of intervals along $\sigma$. 
The $i^{th}$ level of $\mathcal{I}(n)^{\sigma}$ is the collection of its elements of size $i$.
(See Figure \ref{fig_eg} for an example of $\mathcal{I}(6)^{\sigma}$, where $\sigma = 1, 2, 3, 4, 5, 6, 1$.)
For $0 \le i \le n$, let $x_i$ be the number of all subsets of the $i^{th}$ level in $\mathcal{F}$.  
Then 
\[
S = \sum_{\sigma}\left(|\{\emptyset,[n]\} \cap \mathcal{F}|\cdot n + \sum_{i = 1}^{n-1}\binom{n}{i}x_i\right).
\]
Now, we need the following theorem to give an upper bound on $|\{\emptyset,[n]\} \cap \mathcal{F}|\cdot n + \sum_{i = 1}^{n-1}\binom{n}{i}x_i$ for every cyclic permutation $\sigma$.
\begin{theorem}\label{lem_thm}
For every cyclic permutation $\sigma$, we have
\[
|\{\emptyset,[n]\} \cap \mathcal{F}|\cdot n + \sum_{i = 1}^{n-1}\binom{n}{i}x_i \le n\cdot\big(2^{n} - S(n,k,m)\big) + n-1.
\]
\end{theorem}
Supposing that we know Theorem \ref{lem_thm}, then
\[
S=|\mathcal{F}| \cdot n! \le \sum_{\sigma} \Big(n\cdot\big(2^{n} - S(n,k,m)\big) + n-1\Big) = (n-1)! \cdot \Big(n\cdot\big(2^{n} - S(n,k,m)\big) + n-1\Big).
\]
Since $|\mathcal{F}|$ is an integer, we have
\[
|\mathcal{F}| \le \left\lfloor\frac{n\cdot\big(2^{n} - S(n,k,m)\big) + n-1}{n}\right\rfloor = 2^{n} - S(n,k,m),
\]
as desired.
So it is sufficient to prove Theorem \ref{lem_thm}. The following two lemmas are needed to give constraints of $x_0, x_1, \ldots, x_n$. 

\begin{lemma}\label{lem_mid}
For $0 \le i \le n-k+1$, we have
\[
x_i + x_{i+1} + \cdots + x_{i+k-1} \le (k-1)n.
\]
\end{lemma}

\begin{lemma}\label{prop_cons}
If $x_0 = 1$ and $x_1 + x_2 + \cdots +x_{k} = (k-1)n$, then $$x_0 + x_1 + \cdots + x_{k-1} \le (k-1)n - \lfloor n/2\rfloor.$$ 
If $x_n = 1$ and $x_{n-k} + x_{n-k+1} + \cdots +x_{n-1} = (k-1)n$, then $$x_{n-k+1} + x_{n-k+2} + \cdots + x_{n} \le (k-1)n - \lfloor n/2\rfloor.$$
\end{lemma}

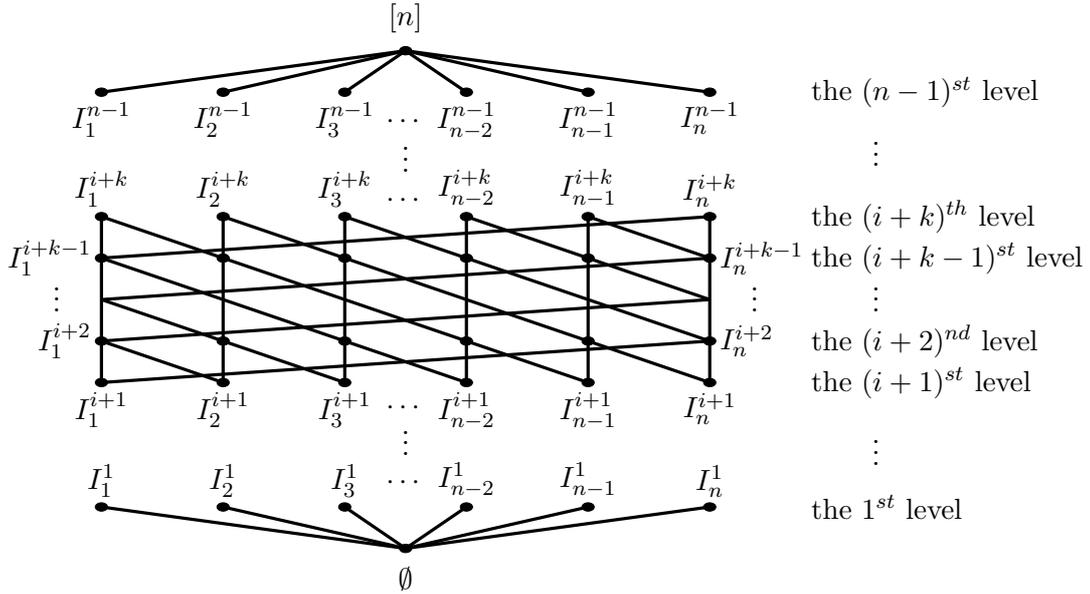
\begin{figure}[htbp]
\begin{center}
\begin{tikzpicture}[xscale=0.8,yscale=0.55]

\draw [fill] (-5,1) circle [radius = 0.1];

\draw [fill] (-3,1) circle [radius = 0.1];

\draw [fill] (-1,1) circle [radius = 0.1];

\draw [fill] (1,1) circle [radius = 0.1];

\draw [fill] (3,1) circle [radius = 0.1];

\draw [fill] (5,1) circle [radius = 0.1];

\draw [fill] (-5,-1) circle [radius = 0.1];

\draw [fill] (-3,-1) circle [radius = 0.1];

\draw [fill] (-1,-1) circle [radius = 0.1];

\draw [fill] (1,-1) circle [radius = 0.1];

\draw [fill] (3,-1) circle [radius = 0.1];

\draw [fill] (5,-1) circle [radius = 0.1];

\draw [fill] (-5,2) circle [radius = 0.1];

\draw [fill] (-3,2) circle [radius = 0.1];

\draw [fill] (-1,2) circle [radius = 0.1];

\draw [fill] (1,2) circle [radius = 0.1];

\draw [fill] (3,2) circle [radius = 0.1];

\draw [fill] (5,2) circle [radius = 0.1];

\draw [fill] (-5,-2) circle [radius = 0.1];

\draw [fill] (-3,-2) circle [radius = 0.1];

\draw [fill] (-1,-2) circle [radius = 0.1];

\draw [fill] (1,-2) circle [radius = 0.1];

\draw [fill] (3,-2) circle [radius = 0.1];

\draw [fill] (5,-2) circle [radius = 0.1];

\node[above] at(-5,2){$I^{i+k}_1$};
\node[above] at(-3,2){$I^{i+k}_2$};
\node[above] at(-1,2){$I^{i+k}_3$};
\node[above] at(0,2){$\cdots$};
\node[above] at(1,2){$I^{i+k}_{n-2}$};
\node[above] at(3,2){$I^{i+k}_{n-1}$};
\node[above] at(5,2){$I^{i+k}_{n}$};
\node[below] at(-5,-2){$I^{i+1}_1$};
\node[below] at(-3,-2){$I^{i+1}_2$};
\node[below] at(-1,-2){$I^{i+1}_3$};
\node[below] at(0,-2.2){$\cdots$};
\node[below] at(1,-2){$I^{i+1}_{n-2}$};
\node[below] at(3,-2){$I^{i+1}_{n-1}$};
\node[below] at(5,-2){$I^{i+1}_{n}$};
\node[left] at(-5,1){$I^{i+k-1}_1$};
\node[left] at(-5.5,0.2){$\vdots$};
\node[left] at(-5,-1){$I^{i+2}_1$};
\node[right] at(5,1){$I^{i+k-1}_n$};
\node[right] at(5.5,0.2){$\vdots$};
\node[right] at(5,-1){$I^{i+2}_n$};

\node[right] at(6.5,2){the $(i+k)^{th}$ level};
\node[right] at(6.5,1){the $(i+k-1)^{st}$ level};
\node[right] at(7.5,0.2){$\vdots$};
\node[right] at(7.5,3.7){$\vdots$};
\node[right] at(7.5,-3.5){$\vdots$};

\node[right] at(6.5,-1){the $(i+2)^{nd}$ level};
\node[right] at(6.5,-2){the $(i+1)^{st}$ level};
\node[right] at(6.5,-5){the $1^{st}$ level};
\node[right] at(6.5,5){the $(n-1)^{st}$ level};
%

%
\draw[][very thick] (-3,1) -- (1,-1); 
\draw[][very thick] (-5,1) -- (-1,-1); 
\draw[][very thick] (-1,1) -- (3,-1); 

\draw[][very thick] (-5,0) -- (-3,-1); 
\draw[][very thick] (1,1) -- (5,-1); 
\draw[][very thick] (3,1) -- (5,0); 

\draw[][very thick] (-5,0) -- (5,1);
\draw[][very thick] (-5,-1) -- (5,0);
\draw[-][very thick] (-5,1) -- (-5,-1); 
\draw[-][very thick] (-3,1) -- (-3,-1); 
\draw[-][very thick] (-1,1) -- (-1,-1); 
\draw[-][very thick] (1,1) -- (1,-1); 
\draw[-][very thick] (3,1) -- (3,-1); 
\draw[-][very thick] (5,1) -- (5,-1); 
%

\draw[-][very thick] (-5,2) -- (-5,1); 
\draw[-][very thick] (-3,2) -- (-3,1); 
\draw[-][very thick] (-1,2) -- (-1,1); 
\draw[-][very thick] (1,2) -- (1,1); 
\draw[-][very thick] (3,2) -- (3,1); 
\draw[-][very thick] (5,2) -- (5,1); 
\draw[-][very thick] (-5,1) -- (5,2); 
\draw[-][very thick] (-5,2) -- (-3,1); 
\draw[-][very thick] (-3,2) -- (-1,1); 
\draw[-][very thick] (-1,2) -- (1,1); 
\draw[-][very thick] (1,2) -- (3,1); 
\draw[-][very thick] (3,2) -- (5,1); 
\draw[-][very thick] (-5,-2) -- (-5,-1); 
\draw[-][very thick] (-3,-2) -- (-3,-1); 
\draw[-][very thick] (-1,-2) -- (-1,-1); 
\draw[-][very thick] (1,-2) -- (1,-1); 
\draw[-][very thick] (3,-2) -- (3,-1); 
\draw[-][very thick] (5,-2) -- (5,-1); 
\draw[-][very thick] (-5,-2) -- (5,-1); 
\draw[-][very thick] (-5,-1) -- (-3,-2); 
\draw[-][very thick] (-3,-1) -- (-1,-2); 
\draw[-][very thick] (-1,-1) -- (1,-2); 
\draw[-][very thick] (1,-1) -- (3,-2); 
\draw[-][very thick] (3,-1) -- (5,-2); 
\draw [fill] (0,6) circle [radius = 0.1];
\draw [fill] (-5,5) circle [radius = 0.1];
\draw [fill] (-3,5) circle [radius = 0.1];
\draw [fill] (-1,5) circle [radius = 0.1];
\draw [fill] (5,5) circle [radius = 0.1];
\draw [fill] (3,5) circle [radius = 0.1];
\draw [fill] (1,5) circle [radius = 0.1];

\draw [fill] (0,-6) circle [radius = 0.1];
\draw [fill] (-5,-5) circle [radius = 0.1];
\draw [fill] (-3,-5) circle [radius = 0.1];
\draw [fill] (-1,-5) circle [radius = 0.1];
\draw [fill] (5,-5) circle [radius = 0.1];
\draw [fill] (3,-5) circle [radius = 0.1];
\draw [fill] (1,-5) circle [radius = 0.1];

\node[above] at(0,6.2){$[n]$};
\node[below] at(0,-6.2){$\emptyset$};

\draw[-][very thick] (0,6) -- (-5,5); 
\draw[-][very thick] (0,6) -- (-3,5); 
\draw[-][very thick] (0,6) -- (-1,5); 
\draw[-][very thick] (0,6) -- (5,5); 
\draw[-][very thick] (0,6) -- (3,5); 
\draw[-][very thick] (0,6) -- (1,5); 

\draw[-][very thick] (0,-6) -- (-5,-5); 
\draw[-][very thick] (0,-6) -- (-3,-5); 
\draw[-][very thick] (0,-6) -- (-1,-5); 
\draw[-][very thick] (0,-6) -- (5,-5); 
\draw[-][very thick] (0,-6) -- (3,-5); 
\draw[-][very thick] (0,-6) -- (1,-5); 

\node[below] at(1,5){$I^{n-1}_{n-2}$};
\node[below] at(3,5){$I^{n-1}_{n-1}$};
\node[below] at(5,5){$I^{n-1}_{n}$};
\node[below] at(0,4.7){$\cdots$};
\node[below] at(-5,5){$I^{n-1}_1$};
\node[below] at(-3,5){$I^{n-1}_2$};
\node[below] at(-1,5){$I^{n-1}_3$};
\node[below] at(0,4.3){$\vdots$};

\node[above] at(1,-5){$I^{1}_{n-2}$};
\node[above] at(3,-5){$I^{1}_{n-1}$};
\node[above] at(5,-5){$I^{1}_{n}$};
\node[below] at(0,-4){$\cdots$};
\node[above] at(-5,-5){$I^{1}_1$};
\node[above] at(-3,-5){$I^{1}_2$};
\node[above] at(-1,-5){$I^{1}_3$};
\node[above] at(0,-4){$\vdots$};

\end{tikzpicture}
\caption{ Vertices $I_t^s$ of $\mathcal{I}(n)^{\sigma}$.}\label{fig_chaindecomp}
\end{center}
\end{figure}

Before starting the proofs of the two lemmas above, we introduce some notations. Let $\sigma = a_1, a_2, \ldots, a_n, a_1 $, where $a_j \in [n]$ for $1 \le j \le n$.  
We denote by $I_{t}^{s}$ the interval $\{a_{t+1}, a_{t+2}, \ldots, a_{t+s}\}$ (with addition taken modulo $n$). 
Clearly, $|I_{t}^{s}| = s$, $I^{s-1}_{t}$, $I^{s-1}_{t+1} \subset I_{t}^{s}$ for $2 \le s \le n-1$, and $I_{t}^{s} \subset I^{s+1}_{t-1}$, $I^{s+1}_{t}$ for $1 \le s \le n-2$ (see Figure \ref{fig_chaindecomp}). 

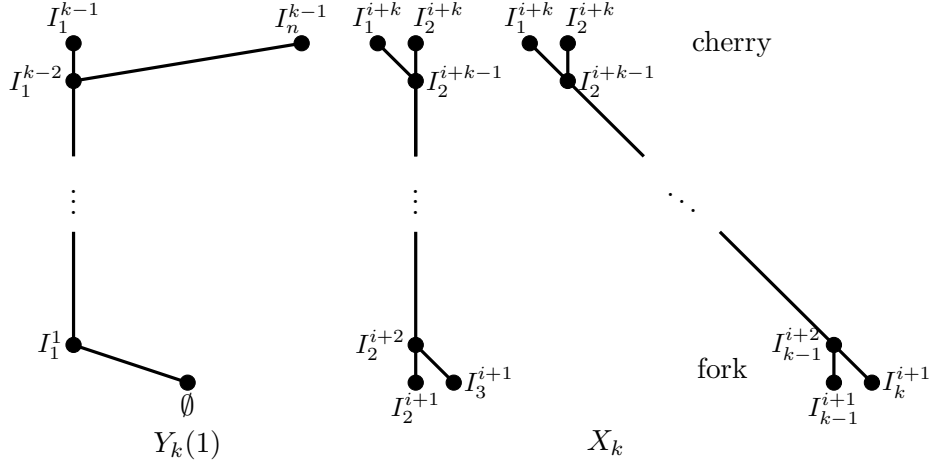
\begin{figure}[htbp]
\begin{center}
\begin{tikzpicture}

\draw [fill] (-7,1) circle [radius = 0.1];
\draw [fill] (-4,1) circle [radius = 0.1];
\draw [fill] (-7,0.5) circle [radius = 0.1];
\draw [fill] (-7,-3) circle [radius = 0.1];
\draw [fill] (-5.5,-3.5) circle [radius = 0.1];

\node[above] at(-7,1){\small{$I_{1}^{k-1}$}};
\node[above] at(-4,1){\small{$I_{n}^{k-1}$}};
\node[left] at(-7,0.5){\small{$I_{1}^{k-2}$}};
\node[] at(-7,-1){$\vdots$};
\node[left] at(-7,-3){\small{$I_{1}^{1}$}};
\node[below] at(-5.5,-3.5){$\emptyset$};
\node[below] at(-5.5,-4){$Y_k(1)$};

\draw[-][very thick] (-7,1) -- (-7,0.5);
\draw[-][very thick] (-7,0.5) -- (-7,-0.5);
\draw[-][very thick] (-7,-1.5) -- (-7,-3);
\draw[-][very thick] (-7,0.5) -- (-4,1);
\draw[-][very thick] (-7,-3) -- (-5.5,-3.5);


\draw [fill] (-3,1) circle [radius = 0.1];
\draw [fill] (-2.5,0.5) circle [radius = 0.1];
\draw [fill] (-2.5,1) circle [radius = 0.1];
\draw [fill] (-2.5,-3.5) circle [radius = 0.1];
\draw [fill] (-2.5,-3) circle [radius = 0.1];
\draw [fill] (-2,-3.5) circle [radius = 0.1];

\draw[-][very thick] (-3,1) -- (-2.5,0.5);
\draw[-][very thick] (-2.5,1) -- (-2.5,0.5);
\draw[-][very thick] (-2.5,0.5) -- (-2.5,-0.5);
\draw[-][very thick] (-2.5,-1.5) -- (-2.5,-3.5);
\draw[-][very thick] (-2.5,-3) -- (-2,-3.5);

\node[] at(-2.5,-1){$\vdots$};
\node[above] at(-3,1){\small{$I_{1}^{i+k}$}};
\node[above] at(-2.2,1){\small{$I_{2}^{i+k}$}};
\node[right] at(-2.5,0.5){\small{$I_{2}^{i+k-1}$}};
\node[below] at(-2.5,-3.5){\small{$I_{2}^{i+1}$}};
\node[left] at(-2.5,-3){\small{$I_{2}^{i+2}$}};
\node[right] at(-2,-3.5){\small{$I_{3}^{i+1}$}};
\node[below] at(0,-4){$X_k$};

\draw [fill] (-1,1) circle [radius = 0.1];
\draw [fill] (-0.5,0.5) circle [radius = 0.1];
\draw [fill] (-0.5,1) circle [radius = 0.1];
\draw [fill] (3,-3.5) circle [radius = 0.1];
\draw [fill] (3,-3) circle [radius = 0.1];
\draw [fill] (3.5,-3.5) circle [radius = 0.1];

\draw[-][very thick] (-1,1) -- (-0.5,0.5);
\draw[-][very thick] (-0.5,1) -- (-0.5,0.5);
\draw[-][very thick] (-2.5,0.5) -- (-2.5,-0.5);
\draw[-][very thick] (-0.5,0.5) -- (0.5,-0.5);
\draw[-][very thick] (1.5,-1.5) -- (3,-3);
\draw[-][very thick] (3,-3) -- (3.5,-3.5);
\draw[-][very thick] (3,-3) -- (3,-3.5);

\node[] at(1,-1){$\ddots$};
\node[above] at(-1,1){\small{$I_{1}^{i+k}$}};
\node[above] at(-0.2,1){\small{$I_{2}^{i+k}$}};
\node[right] at(-0.5,0.5){\small{$I_{2}^{i+k-1}$}};
\node[below] at(3,-3.5){\small{$I_{k-1}^{i+1}$}};
\node[left] at(3,-3){\small{$I_{k-1}^{i+2}$}};
\node[right] at(3.5,-3.5){\small{$I_{k}^{i+1}$}};
\node[left] at(2,-3.3){fork};
\node[right] at(1,1){cherry};


\end{tikzpicture}
\caption{Examples of structures on $k$ consecutive levels.}\label{fig_cherryfork}
\end{center}
\end{figure}

In order to prove the two kinds of constraints of $k$ consecutive  $x_i$'s above, we consider some typical structures in $k$ consecutive  levels of  $\mathcal{I}(n)^{\sigma}$. 
If the $k$ levels are the levels from $0^{th}$ to $(k-1)^{st}$, then every $Y_k$ on consecutive levels must contain $\emptyset$. 
In this case, we consider a special kind of $Y_k$, which we denote by $$Y_k(j) = \{\emptyset, I_j^1, I_j^2, \ldots, I^{k-2}_j, I_j^{k-1}, I_{j-1}^{k-1}\}$$ where $1 \le j \le n$. 
(See Figure \ref{fig_cherryfork} for an example of $Y_k(1)$.) 
If the $k$ levels are the levels from $(i+1)^{st}$ to $(i+k)^{th}$ (the middle part of Figure \ref{fig_chaindecomp}), then we introduce a new kind of structure on $k$ consecutive levels. 
A family of $k+2$ subsets is called $X_k$ if it is $$\{I_{t+1}^{i+1}, I_{t}^{i+1},I_{t}^{i+2}, I_{t}^{i+3}, \ldots, I_{t}^{i+k}, I_{t-1}^{i+k}\} \text{~or~} \{I_{t-1}^{i+1}, I_{t}^{i+1},I_{t-1}^{i+2}, I_{t-2}^{i+3}, \ldots, I_{t-k+1}^{i+k}, I_{t-k+2}^{i+k}\},$$ where $1 \le t \le n.$ 
(See Figure \ref{fig_cherryfork} for two types of examples of $t = 2$ and $t = k$ respectively.) In each $X_k$, we call the $3$ elements in the $(i+k)^{th}$ level and the $(i+k-1)^{st}$ level a cherry, and the $3$ elements in the $(i+1)^{st}$ level and the $(i+2)^{nd}$ level a fork.

\begin{remark}\label{rem_decomp}
$\mathrm{(I)}$ The vertices of the $n$ $Y_k$'s $($$Y_k(1), Y_k(2), \ldots, Y_k(n)$$)$ cover the vertices of the levels from $1^{st}$ to $(k-2)^{nd}$ once. \\
$\mathrm{(II)}$ For $k \ge 4$, the edges of the $2n$ $X_k$'s cover the edges within the $(i+1)^{st}$ level and the $(i+2)^{nd}$ level twice, the edges within the levels from $(i+2)^{nd}$ to $(i+k-1)^{st}$ once, and the edges within the $(i+k-1)^{st}$ level and the $(i+k)^{th}$ level twice.
\end{remark}

Now, we are ready to start the proofs of Lemmas \ref{lem_mid} and \ref{prop_cons}.

\noindent{\bf Proof of Lemma \ref{lem_mid}:~~} 
First, we show $x_0 + x_1 + \cdots + x_{k-1} \le (k-1)n$.  
Since $x_0 \le 1$ and $x_i \le n$ for $1 \le i \le k-1$, we have $x_0 + x_1 + \cdots + x_{k-1} \le (k-1)n + 1$. If $x_0 + x_1 + \cdots + x_{k-1} = (k-1)n + 1$, every subset of the levels from $0^{th}$ to $(k-1)^{st}$ is in $\mathcal{F}$, then one can easily find a copy of $Y_k$ on consecutive levels. 
(See $Y_k(1)$ for an example in Figure \ref{fig_cherryfork}.)
Similarly, we have $x_{n-k+1} + x_{n-k+2} + \cdots + x_{n} \le (k-1)n$.

Now, we prove that
\[
x_{i+1} + x_{i+2} + \cdots + x_{i+k} \le (k-1)n
\]
for $0 \le i \le n-k-1$. 
Recall that there is an edge between two vertices $I, I'$ if and only if $I' \subset I$ and $|I| - |I'| = 1$. Then, we double count the sum $T$ of the weight function $\psi(e)$ over all edges $e = \{I, I'\}$ within the levels from $(i+1)^{st}$ to $(i+k)^{th}$ (the middle part of Figure \ref{fig_chaindecomp}), where
\begin{align*}
\psi(e) =
\begin{cases}
0, &\mbox{if $I, I' \in \mathcal{F}$;}\\
0, &\mbox{if $I, I' \notin \mathcal{F}$;}\\
0, &\mbox{$I\in \mathcal{F}$, $I' \notin \mathcal{F}$ and $|I| = i+k$; }\\
2, &\mbox{$I\notin \mathcal{F}$, $I' \in \mathcal{F}$ and $|I| = i+k$; }\\
0, &\mbox{$I\notin \mathcal{F}$, $I' \in \mathcal{F}$ and $|I| = i+2$; }\\
2, &\mbox{$I\in \mathcal{F}$, $I' \notin \mathcal{F}$ and $|I| = i+2$; }\\
1, &\mbox{otherwise.}
\end{cases}
\end{align*} 

On the one hand, recall that every vertex in the $s^{th}$ level of $\mathcal{I}(n)^{\sigma}$ have two neighbors in the $(s+1)^{st}$ level and two neighbors in the $(s-1)^{st}$ level, and note that only the case when
one of the two vertices in the edge is not in $\mathcal{F}$ counts nonzero. 
So
\begin{align}\label{prop_most}
    T=&\sum_{I \notin \mathcal{F}} \psi(e) + \sum_{I' \notin \mathcal{F}} \psi(e)\notag \\
\le&\Big(2\cdot2(n - x_{i+k}) +1\cdot \sum_{j = i + 3}^{i+k-1}2(n-x_j)\Big)
     + \Big(2\cdot2(n - x_{i+1}) +1\cdot\sum_{l = i+2}^{i+k-2} 2(n - x_{l})\Big)\notag\\ 
 =&4\Big(\sum_{h = i+1}^{i+k}(n - x_{h})\Big) - 2(n - x_{i+2}) - 2(n -x_{i+k-1})\notag.
\end{align}

On the other hand, we will show  $T \ge 2(x_{i+2}+x_{i+k-1})$.
When $k = 3$, since every vertex of the $(i+2)^{nd}$ level in $\mathcal{F}$ has at least $2$ neighbors not in $\mathcal{F}$, it follows that $T \ge 2\cdot 2\cdot x_{i+2}.$
If $k \ge 4$, we consider the $2n$ $X_k$'s, and change the weight function of an edge to $1$ if it originally counts $2$. Then by (II) of Remark \ref{rem_decomp}, $T$ is equal to the weighted sum over all edges in these $2n$ $X_k$'s.
We will prove below the claim that in each $X_k$, the sum of the changed weight function over its edges  is at least $|\{I \in \mathcal{F}\cap X_k \mid \text{$|I| = i+2$ or  $|I| = i+ k- 1$}\}|.$
Then if we summarize it over all $2n$ $X_k$'s, $T$ is at least $2(x_{i+2}+x_{i+k-1})$.  

\begin{figure}[htbp]
\begin{center}
\begin{tikzpicture}


\draw [fill] (-7,1) circle [radius = 0.1];
\draw [fill] (-5,1) circle [radius = 0.1];
\draw [fill] (-3,1) circle [radius = 0.1];
\draw [] (-1,1) circle [radius = 0.1];
\draw [] (1,1) circle [radius = 0.1];
\draw [] (3,1) circle [radius = 0.1];
\draw [dashed] (5,1) circle [radius = 0.1];
\draw [dashed] (7,1) circle [radius = 0.1];
\draw [fill] (-6,0) circle [radius = 0.1];
\draw [fill] (-2,0) circle [radius = 0.1];
\draw [fill] (2,0) circle [radius = 0.1];
\draw [] (6,0) circle [radius = 0.1];

\draw [fill] (-7,-3) circle [radius = 0.1];
\draw [fill] (-5,-3) circle [radius = 0.1];
\draw [fill] (-3,-3) circle [radius = 0.1];
\draw [] (-1,-3) circle [radius = 0.1];
\draw [] (1,-3) circle [radius = 0.1];
\draw [] (3,-3) circle [radius = 0.1];
\draw [dashed] (5,-3) circle [radius = 0.1];
\draw [dashed] (7,-3) circle [radius = 0.1];
\draw [fill] (-6,-2) circle [radius = 0.1];
\draw [fill] (-2,-2) circle [radius = 0.1];
\draw [fill] (2,-2) circle [radius = 0.1];
\draw [] (6,-2) circle [radius = 0.1];

\node[below] at(-6,-0.2){Type 1};
\node[below] at(-2,-0.2){Type 2};
\node[below] at(2,-0.2){Type 3};
\node[below] at(6,-0.2){Type 4};

\node[below] at(-6,-3.2){Type 5};
\node[below] at(-2,-3.2){Type 6};
\node[below] at(2,-3.2){Type 7};
\node[below] at(6,-3.2){Type 8};

\draw[-][very thick] (-7,1) -- (-6,0);
\draw[-][very thick] (-5,1) -- (-6,0);
\draw[-][very thick] (-3,1) -- (-2,0);
\draw[-][very thick] (-1,0.9) -- (-2,0);
\draw[-][very thick] (1,0.9) -- (2,0);
\draw[-][very thick] (3,0.9) -- (2,0);
\draw[-][very thick] (5,0.9) -- (6,0.1);
\draw[-][very thick] (7,0.9) -- (6,0.1);

\draw[-][very thick] (-7,-3) -- (-6,-2);
\draw[-][very thick] (-5,-3) -- (-6,-2);
\draw[-][very thick] (-3,-3) -- (-2,-2);
\draw[-][very thick] (-1,-2.9) -- (-2,-2);
\draw[-][very thick] (1,-2.9) -- (2,-2);
\draw[-][very thick] (3,-2.9) -- (2,-2);
\draw[-][very thick] (5,-2.9) -- (6,-2.1);
\draw[-][very thick] (7,-2.9) -- (6,-2.1);

\end{tikzpicture}
\caption{Types $1-8$: the vertices $I$ in $\mathcal{F}$ (not in $\mathcal{F}$ or not clear) are denoted by solid (hollow or dashed) vertices, respectively.}\label{fig_types}
\end{center}
\end{figure}
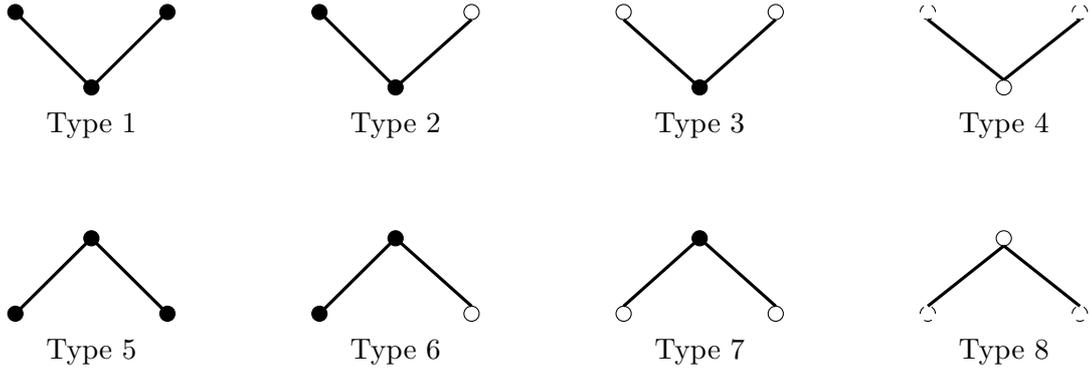

Now, we divide the cherries of all $2n$ $X_k$'s into $4$ types (see Figure \ref{fig_types}), according to whether its $3$ elements are in $\mathcal{F}$ or not. 
We call a cherry Type $1$ $(2$ or $3)$ if the middle element and both (one or none) of its neighbors are in $\mathcal{F}$. The rest cases are Type $4$, namely that the middle element is not in $\mathcal{F}$. Note that in Type $4$, we do not distinguish the cases if two neighbors of the middle element are in $\mathcal{F}$ or not. Similarly, we divide all forks  into Types $5, 6, 7$ and $8$ (see Figure \ref{fig_types}).

To prove the claim, we distinguish several cases by the types defined above. 
Note that in each $X_k$, in order to avoid copies of $Y_k$ and $Y'_k$ in $\mathcal{F}$, either at least $2$ of the $4$ vertices in the $(i+1)^{st}$ level and the $(i+k)^{th}$ level are not in $\mathcal{F}$, or at least one of the $k-2$ vertices in the levels from $(i+2)^{nd}$ to $(i+k-1)^{st}$ is not in $\mathcal{F}$.
If the cherry and the fork in some $X_k$ are Types $1$ and $5$ respectively, then we must have some vertex $I$ not in $\mathcal{F}$ in the chain connecting the cherry and the fork in this $X_k$. 
Then along this chain in both directions from $I$, take the first vertices in $\mathcal{F}$  respectively.  Thus, we can find at least  two edges with weight $1$. 
Note that such vertices exist since in Types $1$ and $5$, the vertices in  the $(i+2)^{nd}$ level and the $(i+ k -1)^{st}$ level are in $\mathcal{F}$. 
The same argument works for the pairs Types $1$ and $6$, and Types $2$ and $5$. 
If they are Types $1$ and $7$, then we already have two desired edges in the fork. 
If they are Types $1$ and $8$, then we find the first vertex in $\mathcal{F}$ along the chain from the fork to the cherry, this vertex and its neighbor below in the chain will give the edge we want.
The same argument works for Types $4$ and $5$.
In the rest cases, it is easy to find the desired edges in the cherry and the fork. (If $k =4$, Types $1$ and $5$ ($1$ and $6$, or $2$ and $5$) cannot appear in any $X_k$, since they will form a $Y_k$ or a $Y'_k$ in $\mathcal{F}$.)

Therefore, we have 
\[
2(x_{i+2} + x_{i+k-1}) \le  4\Big(\sum_{h = i+1}^{i+k}(n - x_{h}) \Big) - 2(n - x_{i+2}) - 2(n -x_{i+k-1}).
\]
That is,
\[
x_{i+1} + x_{i+2} + \cdots + x_{i+k} \le (k-1)n,
\]
as required. $\hfill \Box$

\noindent{\bf Proof of Lemma \ref{prop_cons}:~~} By symmetry, it is enough to prove the first part of the lemma. 
That is, we need to prove 
\[
x_0 + x_1 + \cdots + x_{k-1} \le (k-1)n - \lfloor n/2 \rfloor = (k -2)n + \lceil n/2\rceil.
\]
We distinguish three cases to prove it.
\begin{case}
 $x_{k-1} \le \lfloor n/2\rfloor$ when $n$ is odd or $x_{k-1} \le  (n/2)-1$ when $n$ is even.
\end{case}

Note that $x_i \le n$ for $i = 1, 2, \ldots, k-2$, then by the assumption $x_0 = 1$, we have
\[
x_0 + x_1 + \cdots + x_{k-1} \le 1  + (k-2)n +  \lfloor n/2\rfloor  \le (k -2)n + \lceil n/2\rceil,
\]
if $n$ is odd; and
\[
x_0 + x_1 + \cdots + x_{k-1} \le 1  + (k-2)n +   (n/2) - 1 \le (k -2)n + (n/2),
\]
if $n$ is even.

\begin{case}
$x_{k-1} =  n/2$ when $n$ is even.
\end{case}
Suppose that the result is not true. That is, $$x_0 + x_1 + \cdots + x_{k-1} \ge (k -2)n+(n/2)+1.$$ Note that $x_{k-1} =  n/2$ in this case and $x_i \le n$ for $i = 1, 2, \ldots, k-2$. Then  by the assumption $x_0 = 1$, we have  $x_1 = x_2 = \cdots = x_{k-2} = n$. Therefore, by the assumption $x_1 + x_2 + \cdots +x_{k} = (k-1)n$, we have $x_{k} =  n/2$. 
Since $x_{k-1} =  n/2$ and $x_{k} =  n/2$, we can find $2$ vertices such that they are in the $k^{th}$ level and the $(k-1)^{st}$ level respectively, they are in $\mathcal{F}$, and they have inclusion relation. Then a copy of $Y'_k$ on the levels from $1^{st}$ to $k^{th}$ can be found, a contradiction. 
\begin{case}
$x_{k-1} \ge \lfloor n/2\rfloor + 1$. 
\end{case}

\begin{figure}[htbp]
\begin{center}
\begin{tikzpicture}
\draw (0,0) circle (1); 
\draw [fill] (-0.866,-0.5) circle [radius = 0.1];
\draw [fill] (0.866,-0.5) circle [radius = 0.1];
\draw [] (0.866,0.5) circle [radius = 0.1];
\draw [fill] (0,-1) circle [radius = 0.1];
\draw [fill] (0,1) circle [radius = 0.1];

\draw [fill] (1,0) circle [radius = 0.1];
\draw [] (-0.5,-0.866) circle [radius = 0.1];
\draw [fill] (-0.5,0.866) circle [radius = 0.1];
\draw [] (0.5,-0.866) circle [radius = 0.1];
\draw [] (0.5,0.866) circle [radius = 0.1];
\node[left] at(-0.866,0.5){$\cdot$};
\node[left] at(-0.97,0.25){$\cdot$};
\node[left] at(-0.99,0){$\cdot$};
\node[right] at(0.5,1){$I^{k-1}_{1}$};
\node[above] at(0,1){$I^{k-1}_{n}$};
\node[left] at(-0.5,1){$I^{k-1}_{n-1}$};
\end{tikzpicture}
\caption{The cycle of the vertices in the $(k-1)^{st}$ level.}\label{fig_cycle}
\end{center}
\end{figure}
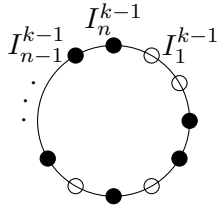

In this case, we count the number of pairs of vertices of the $(k-1)^{st}$ level in $\mathcal{F}$ which are neighbors (having a common neighbor vertex in the $(k-2)^{nd}$ level). 
Now, we view the $n$ vertices in the $(k-1)^{st}$ level  around a cycle (see Figure \ref{fig_cycle}), and call a set of vertices  a part if these vertices are in $\mathcal{F}$ and consecutive in the cycle.
Note that the number of the vertices in the $(k-1)^{st}$ level which are not in $\mathcal{F}$ is $n-x_{k-1}$, and these vertices cut the cycle into at most $n-x_{k-1}$ parts. 
In every part $P$, the number of the pairs we want is equal to  $|P|-1.$
So in total, we have at least $$x_{k-1} - (n-x_{k-1}) = 2x_{k-1} - n$$ pairs of vertices in the  $(k-1)^{st}$ level, which may form a $Y_k(j)$ in $\mathcal{F}$ (see Figure \ref{fig_cherryfork}). Then in each of such $2x_{k-1} - n$ candidates ($Y_k(j)$'s), at least one vertex in the levels from $1^{st}$ to $(k-2)^{nd}$ is not in $\mathcal{F}$. So by (I) of Remark \ref{rem_decomp}, we have
\[
\sum_{i =0}^{k-1}x_i \le 1+  (k-2)n - (2x_{k-1} - n)+ x_{k-1} = 1+ (k-1)n - x_{k-1} \le (k-1)n - \lfloor n/2 \rfloor,
\]
since $x_{k-1} \ge \lfloor n/2\rfloor + 1$. This completes the proof of this case and the lemma. $\hfill \Box$

Before starting the proof of Theorem \ref{lem_thm}, we need the following notations and lemmas for helping us to use the constraints above. Recall the definition of $S(n,k,r)$ and $m=\lceil(n-k)/2\rceil$. Now let $z$ be an integer such that $z \equiv r$ (mod $k$). 
Then we denote 
\begin{align*}
S(n,k,r\mid z) = 
 \sum_{\substack{i=0\\ i \equiv r \text{~mod~} k\\i \le z}}^{n}\binom{n}{i};
\end{align*}
\begin{align*}
w_{i}=
\begin{cases}
S(n,k,i\mid i) - S(n,k,i-1\mid i-1), &\mbox{if $i \le m$;} \\ 
S(n,k,n-i-k+1\mid n-i-k+1) - S(n,k,n-i-k\mid n-i-k), &\mbox{if $m + 1\le i$.}\\
\end{cases}
\end{align*}

\begin{lemma}\label{lem:lem3} Let $t \le n$ be an integer. We have the following equalities. \\
$\mathrm{(I)}$ If $t < 0$, $S(n,k,t \mid t) = 0$,  and if $0 \le t \le k-1$, $S(n,k,t\mid t) = \binom{n}{t}$.\\
$\mathrm{(II)}$ If $0 \le t \le n$, $S(n,k,t\mid t) = S(n,k,t-k\mid t-k) + \binom{n}{t}$.\\
$\mathrm{(III)}$ If $0 \le t \le n$,
$S(n,k,t\mid t) + S(n,k,n-k-t\mid n-k-t) = S(n,k,t)$.\\
$\mathrm{(IV)}$ If $0 \le t \le n$,
$S(n,k,t-k\mid t-k) + S(n,k,n - k - t\mid n-k-t) = S(n,k,t) - \binom{n}{t}$.
\end{lemma}
\noindent{\bf Proof:~~} $\mathrm{(I)}$ follows from the definition of $S(n,k,t\mid t)$. $\mathrm{(II)}$ follows from $\mathrm{(I)}$ if $0 \le t \le k-1$, and the definition of $S(n,k,t\mid t)$ if $k \le t \le n$. $\mathrm{(IV)}$ follows from $\mathrm{(II)}$ and $\mathrm{(III)}$. Now, we distinguish two cases to prove $\mathrm{(III)}$. 

If $t > n-k$, then we have $t + k > n$ and  $n-k-t <0$. Hence, $S(n,k,t|t) = S(n,k,t)$ and $S(n,k,n - k - t|n-k-t) = 0$ by $\mathrm{(I)}$, and this completes the proof of this case.

If $0 \le t \le n-k$, let $0 \le r \le k-1$ be the remainder of $n-k-t$ divided by $k$. We have 
\begin{align*}
 &  S(n,k,t|t) + S(n,k,n - k - t|n-k-t)  \\
=&  S(n,k,t|t) + \binom{n}{r} + \binom{n}{r+k} + \cdots + \binom{n}{n-k-t}\\
=&  S(n,k,t|t) + \binom{n}{n-r} + \binom{n}{n-r-k} + \cdots + \binom{n}{k+t}\\ 
=&  S(n,k,t|t) + \binom{n}{t+k} + \cdots + \binom{n}{n-r-k} + \binom{n}{n-r}= S(n,k,t),
\end{align*}
since $r < k$. \noindent $\hfill \Box$

\setcounter{case}{0}

\begin{lemma}\label{lem_weight}
Let $s$ be an integer.  \\
$\mathrm{(I)}$ $w_s =0$ if $s <0 $ or $s >n-k+1$, and  $w_s > 0$ if $0 \le s \le n-k+1$.

\noindent$\mathrm{(II)}$ If $0 \le s \le m$ or $m+k \le s \le n$, then
\[
\sum_{i = s-k+1}^{s}w_i = \binom{n}{s}.
\]

\noindent$\mathrm{(III)}$ If $m < s < m+k$,
\begin{align*}
\sum_{i = s-k+1}^{s}w_i= S(n,k,m) - S(n,k,s) + \binom{n}{s}.
\end{align*}

\noindent$\mathrm{(IV)}$ $\sum_{i=0}^{n-k+1}w_i=S(n,k,m)$.
\end{lemma}
\noindent{\bf Proof:~~} $\mathrm{(I)}$ It follows from the definition of $w_s$ and $\mathrm{(I)}$ of Lemma \ref{lem:lem3}.

\noindent$\mathrm{(II)}$ By $\mathrm{(II)}$ of Lemma \ref{lem:lem3}, if $0 \le s \le m$, it follows from $\mathrm{(I)}$ that
\[
\sum_{i = s-k+1}^{s}w_i = \sum_{i = s-k+1}^{s}\big(S(n,k,i|i) - S(n,k,i-1|i-1)\big) = S(n,k,s|s) -S(n,k,s-k|s-k) = \binom{n}{s},
\]
and if $m+k \le s \le n$, we have $s - k + 1 \ge m+1$, and so
\begin{align*}
\sum_{i = s-k+1}^{s}w_i = &\sum_{i = s-k+1}^{s}\big(S(n,k,n-k-i+1|n-k-i+1) - S(n,k,n-k-i|n-k-i)\big) \\= & S(n,k,n-s|n-s) - S(n,k,n-s-k|n-s-k) = \binom{n}{n-s}=\binom{n}{s}.
\end{align*}

\noindent$\mathrm{(III)}$ If $m < s < m+k$, then $s - k + 1 \le m$ and $s \ge m+1$. So by $\mathrm{(III)}$ and $\mathrm{(IV)}$ of Lemma \ref{lem:lem3},            we have 
\begin{align*}
  \sum_{i = s-k+1}^{s}w_i  =&  \sum_{i = s-k+1}^{m}\big(S(n,k,i|i) - S(n,k,i-1|i-1)\big) \\&+ \sum_{j = m+1}^{s}\big(S(n,k,n-k-j+1|n-k-j+1) - S(n,k,n-k-j|n-k-j)\big)\\=& S(n,k,m|m) -S(n,k,s-k|s-k) \\&+ S(n,k,n-k-m|n-k-m) - S(n,k,n-k-s|n-k-s)\\
                         = &\big(S(n,k,m|m) + S(n,k,n-k-m|n-k-m)\big) \\ 
                             &- \big(S(n,k,s - k|s-k) + S(n,k,n-k-s|n-k-s)\big)\\
= &S(n,k,m) - \left(S(n,k,s) - \binom{n}{s}\right)= S(n,k,m) - S(n,k,s) + \binom{n}{s}.
\end{align*}

\noindent$\mathrm{(IV)}$ 
\[
\sum_{i=0}^{n-k+1}w_i = \sum_{\substack{i=0\\ i \equiv m \text{~mod~} k}}^{n}\left(\sum_{j = i-k+1}^{i}w_j\right)=S(n,k,m),
\]
by $\mathrm{(I)}$ and $\mathrm{(II)}$.$\hfill \Box$

\section{Proof of Theorem \ref{lem_thm}}

First, we consider a linear programming problem: maximize  $\sum_{i = 0}^{n}\binom{n}{i}x_i$ subject to the constraints $x_i + x_{i+1} + \cdots + x_{i+k-1} \le (k-1)n$ 
for $i =0, 1 \ldots, n-k+1$ (Lemma \ref{lem_mid}), $0 \le x_0, x_n \le 1$, and $0 \le x_i \le n$ for all $i = 1, \ldots, n-1$.

We assign the weight $w_i$ (see the definition and properties of $w_i$ in Section $2$) to the constraint $x_i +  x_{i+1} + \cdots + x_{i+k-1} \le (k-1)n$ for every $0 \le i \le n-k+1$. By $\mathrm{(I)}$ of Lemma \ref{lem_weight}, $w_i > 0$ for $0 \le i \le n-k+1$, so we have 
\begin{align}\label{eq_weight}
\sum_{i=0}^{n-k+1} w_i(x_i +  x_{i+1} + \cdots + x_{i+k-1})  \le (k-1)n\cdot\sum_{i =0}^{n-k+1}  w_{i}.
\end{align}
Then by $\mathrm{(I)}$, $\mathrm{(II)}$ and $\mathrm{(III)}$ of Lemma \ref{lem_weight}, the LHS of  (\ref{eq_weight}) is equal to
\begin{align*}
   \sum_{i = 0}^{n}\left(\left(\sum_{j=i-k+1}^{i}w_j\right)x_i\right)
=  \sum_{i=0}^{n} \binom{n}{i}x_i + \sum_{j=m + 1}^{m+k-1}\left(S(n,k,m) - S(n,k,j)\right)x_j.
\end{align*}
On the other hand, by $\mathrm{(IV)}$ of Lemma \ref{lem_weight}, the RHS of (\ref{eq_weight}) is equal to $(k-1)n\cdot S(n,k,m)$. Hence, by (\ref{eq_weight}), we have
\begin{align}\label{eq_final}
\sum_{i=0}^{n} \binom{n}{i}x_i + \sum_{j=m + 1}^{m+k-1}\left(S(n,k,m) - S(n,k,j)\right)x_j \le (k-1)n\cdot S(n,k,m).
\end{align}
Note that $S(n,k,m) = \min_{r: \ 0\le r \le k-1}S(n,k,r)$ by Remark \ref{rem_lm}, and $0 \le x_j \le n$ for $m+1 \le j \le m+k-1$. So we have
\begin{align}\label{eq_con}
\sum_{j=m + 1}^{m+k-1}\left(S(n,k,j) - S(n,k,m)\right)x_j \le \sum_{j=m + 1}^{m+k-1}\left(S(n,k,j) - S(n,k,m)\right)n.
\end{align}
Therefore, combining (\ref{eq_final}) and (\ref{eq_con}), we have
\begin{align}\label{eq_sum}
\sum_{i=0}^{n} \binom{n}{i}x_i &\le (k-1)n\cdot S(n,k,m)+\sum_{j=m + 1}^{m+k-1}\left(S(n,k,j) - S(n,k,m)\right)n\notag\\
                               & =  \sum_{j=m + 1}^{m+k-1}\left(S(n,k,j) - S(n,k,m) + S(n,k,m)\right)n \notag\\ &= n\cdot \sum_{j=m+1}^{m+k-1}S(n,k,j)= n\cdot \left(2^n - S(n,k,m)\right).
\end{align}

Now, we use (\ref{eq_sum}) to prove Theorem \ref{lem_thm}.
The first case $|\mathcal{F} \cap \{\emptyset, [n]\}| = 0$ follows from (\ref{eq_sum}).
Then if $|\mathcal{F} \cap \{\emptyset, [n]\}| = 1$, we may suppose that $\emptyset \in \mathcal{F}$ and $[n] \notin \mathcal{F}$, and so $x_0 = 1$ and $x_n =0$. Thus, it follows that
\[
       n\cdot x_0 + \sum_{i = 1}^{n-1}\binom{n}{i}x_i + n\cdot x_n
   =  \sum_{i = 0}^{n}\binom{n}{i}x_i + n -1
  \le n \cdot \big(2^n - S(n,k,m)\big) + n-1,
\]
as required. 

If $\emptyset \in \mathcal{F}$ and $[n] \in \mathcal{F}$, we have  $x_0 = 1$ and $x_n = 1$.
By Lemma \ref{lem_mid}, $x_1 + x_2 + \cdots + x_k \le (k-1)n$ and $x_{n-k} + x_{n-k+1} + \cdots + x_{n-1} \le (k-1)n$. According to these two constraints, we distinguish two subcases.

If $x_1 + x_2 + \cdots + x_k \le (k-1)n - 1$ or $x_{n-k} + x_{n-k+1} + \cdots + x_{n-1} \le (k-1)n - 1$, we may assume $x_1 + x_2 + \cdots + x_k \le (k-1)n - 1$. Then in this case, (\ref{eq_weight}) should be 
\begin{align*}
\sum_{i=0}^{n-k+1} w_i(x_i +  x_{i+1} + \cdots + x_{i+k-1})  \le \left((k-1)n\cdot\sum_{i =0}^{n-k+1}  w_{i}\right) - w_1.
\end{align*}
Since $w_1 = n-1$, (\ref{eq_sum}) should be $\sum_{i=0}^{n} \binom{n}{i}x_i \le n\cdot \left(2^n - S(n,k,m)\right) - (n-1)$, and so
\begin{align*}
    &n(x_0 + x_n) + \sum_{i = 1}^{n-1}\binom{n}{i}x_i
=   (n -1)(x_0+x_n) + \sum_{i = 0}^{n}\binom{n}{i}x_i\\
 \le &2(n-1) + n \cdot \big(2^n - S(n,k,m)\big) - (n-1)\\ 
   = &n \cdot \big(2^n - S(n,k,m)\big) + (n-1).
\end{align*}

If $x_1 + x_2 + \cdots + x_k = (k-1)n$ and $x_{n-k} + x_{n-k+1} + \cdots + x_{n-1} = (k-1)n$, then by Lemma \ref{prop_cons}, $x_0 + x_1 + \cdots + x_{k-1} \le (k-1)n - \lfloor\frac{n}{2}\rfloor$
and $x_{n-k+1} + x_{n-k+2} + \cdots + x_{n} \le (k-1)n - \lfloor\frac{n}{2}\rfloor.$
Similarly, one can modify (\ref{eq_weight}) and (\ref{eq_sum}) to show that
\begin{align*}
n\cdot (x_0+x_n) + \sum_{i = 1}^{n-1}\binom{n}{i}x_i &\le
2(n-1) + n\cdot \big(2^n - S(n,k,m)\big) - 2\lfloor\frac{n}{2}\rfloor\\
& \le (n-1) + n\cdot \big(2^n - S(n,k,m)\big)
\end{align*}
by $w_0 = w_{n-k+1} =1$. This completes the proofs of  Theorems \ref{lem_thm} and  \ref{thm}.$\hfill \Box$

\end{document}